
%
%
%
%
%
\input amstex
%
%
\def\cstar{$C^*$-algebra}

\def\deg{\roman{deg }}
\def\dim{\roman{dim}}

%
\documentstyle{amsppt}
\loadbold
\magnification=\magstep 1
\rightheadtext{infinite dimensional
numerical linear algebra}
\topmatter
\title
The role of $C^*$-algebras in infinite
dimensional numerical linear algebra
\endtitle
\author William Arveson
\endauthor
\affil
Department of Mathematics\\
University of California\\
Berkeley, CA 94720 USA
\endaffil
\email
arveson\@math.berkeley.edu
\endemail
\date
6 June, 1993
\enddate
\rightheadtext{$C^*$-algebras in numerical linear algebra}
\thanks
This research was supported in part by
NSF grant DMS89-12362
\endthanks
\abstract
This paper deals with mathematical issues relating to
the computation of spectra of self adjoint
operators on Hilbert spaces.  We describe
a general method for approximating the spectrum
of an operator $A$ using the eigenvalues
of large finite dimensional
truncations of $A$.  The results of several papers
are summarized which imply that the method is
effective in most cases of interest.
Special attention is paid to
the Schr\"odinger operators of one-dimensional
quantum systems.

We believe that these results serve to make
a broader point, namely that numerical problems involving
infinite dimensional operators
require a reformulation in terms of
\cstar s.  Indeed, it is only when the given operator
$A$ is viewed as an element
of an appropriate \cstar\ $\Cal A$
that one can see the precise
nature of the limit of the finite dimensional
eigenvalue distributions: the
limit is associated with a tracial state on $\Cal A$.
For example, in the case where $A$ is the discretized
Schr\"odinger operator associated with a
one-dimensional quantum system, $\Cal A$ is a simple
\cstar\ having a unique tracial state.  In these cases
there is a precise asymptotic result.
\endabstract

\endtopmatter
\vfill\eject
\pagebreak
\document
\subheading{1. Introduction}
We  discuss methods
for computing the spectrum $\sigma(A)$, and especially the
essential spectrum $\sigma_e(A)$, of bounded self-adjoint
operators $A$ on separable Hilbert spaces.

We mean to take the term `compute' seriously.
Most operators
that arise in practice are not
presented in a representation
in which they are diagonalized, and
it is often very
hard to locate even a single
point in the spectrum
of the operator \cite{6}, \cite{7}, \cite{8}, \cite{12},
 \cite{13}, \cite{15}.  Some typical examples
will be discussed in section 3.  Thus, one often has
to settle for numerical approximations to  $\sigma(A)$ or
$\sigma_e(A)$, and this raises the question of how to
implement the methods of finite dimensional
numerical linear algebra to compute the spectra of infinite
dimensional operators.  Unfortunately, there is a dearth of
literature on this basic problem and, so far as we have been
able to tell, there are no {\it proven} techniques.

In this paper we establish an effective method
for approaching such
problems; we discuss issues
associated with operator theory and operator algebras,
but not issues
belonging properly to numerical analysis.
Thus, we address the question of whether or
not certain finite dimensional approximations converge
to the correct limit, but
we do not address questions relating
to how fast they converge nor how this method might be
implemented algorithmically.  Nevertheless, it may
be appropriate to point out that these
methods have been effectively implemented in a Macintosh
program \cite{2} which is available from the author.

The following section contains a general method for
computing the essential spectrum of a self adjoint operator.
The method depends on choosing an orthonormal
basis and approximating the
operator with $n\times n$ matrices
obtained as sections of the
infinite matrix of the operator
with respect to this basis.
Not every orthonormal basis
is appropriate, and we develop criteria
which show how the basis must be chosen.
This involves an
abstraction of the classical notion
of band-limited matrix.

Our method is
effective for computing the
essential spectrum.  When
there is a difference between
the spectrum and essential
spectrum, our results give no
information about points
in the difference $\sigma(A)\setminus\sigma_e(A)$.
However, since $\sigma(A)\setminus\sigma_e(A)$
consists
merely of isolated eigenvalues
of finite multiplicity,
this difference is usually
insignificant and
often empty (as it is for the main
examples below).

In section 3 we describe the source
of our principal
examples of self adjoint operators.
We show how one
should `discretize' the canonical
commutation relations
so as to preserve the uncertainty
principle, and relate
the resulting Schr\"odinger operators
to the tridiagonal
operators considered in section 4.
The one-parameter
family of \cstar s generated
by the discretized canonical commutation
relations turns out to be the
one-parameter family of noncommutative
spheres of Bratteli, Elliott,
Evans and Kishimoto
\cite{9}, \cite{10}, \cite{11}.
In our case the parameter is related to the
numerical step size.

In section 4 we clarify the role of
\cstar s in numerical
problems of this kind.  In particular,
the \cstar\ s
associated with a broad class of
tridiagonal operators are
simple \cstar\ s having a unique
tracial state.  The tracial
state plays an essential role in
describing the limit of
the eigenvalue distributions of
the approximating sequence
of finite dimensional truncations
of the basic operator.

While there are some new results below,
what follows is primarily an exposition of
the results of several recent papers
\cite{1}, \cite{2}, \cite{3}, \cite{5}.

\subheading{2. Filtrations and degree}
Let $H$ be a Hilbert space.  A
{\sl filtration} of $H$
is an increasing sequence
$$
\Cal F = \{H_1\subseteq H_2\subseteq \dots \}
$$
of finite dimensional subspaces
$H_n$ of $H$ with the
property that
$$
\lim_{n\to \infty}\dim H_n = \infty.
$$
The filtration $\Cal F$ is called
{\sl proper} or
{\sl improper} according as the union
$\cup H_n$ is
dense in $H$ or is not dense in $H$.
In general, we
will write $P_n$ and $P_+$ for
the projections onto
the subspaces $H_n$ and
$\overline{\cup H_n}$
respectively.

The simplest
filtrations are associated with
orthonormal sets
in $H$.  For example, if
$\{e_n: n=1,2,\dots\}$ is
an orthonormal set then
$$
H_n = [e_1, e_2, \dots, e_n], \qquad n=1,2,\dots \tag{2.1}
$$
defines a filtration of $H$ which is proper iff
$\{e_n\}$ is an orthonormal basis.
If, on the other hand, we are given
a bilateral orthonormal basis
$\{e_n: n=0, \pm 1,\pm 2,\dots\}$
for $H$, then
$$
H_n = [e_n, e_{-n+1},\dots, e_n]
$$
defines a proper filtration in which the
dimensions increase
in jumps of $2$.  Moreover, there is a natural
{\sl improper} filtration associated
with such a bilateral basis,
namely $\Cal F^+=\{H_1^+\subseteq H_2^+\subseteq\dots\}$
where
$$
H_n^+ = [e_0,e_1,\dots,e_n],\qquad n=1,2,\dots.
$$
The spaces $H_n^+$ span a proper subspace
$H_+$ of $H$.
This example of an improper filtration
is associated with
`unilateral' sections.  That is to say,
the matrix of any
operator $A\in\Cal B(H)$
relative to the basis
$\{e_n: n\in \Bbb Z\}$ is a doubly
infinite matrix $(a_{ij})$,
whereas the matrix of the compression
$$
P_+A\restriction_{H_+}
$$
is a singly infinite submatrix of
$(a_{ij})$.  While these
are the main examples of filtrations
for our purposes here,
filtrations in general are
allowed to have very
irregular jumps in dimension.

Given any filtration
$\Cal F = \{H_n: n\geq1\}$ of $H$,
we introduce the notion of degree
of an operator
(relative to $\Cal F$) as follows.

\proclaim{Definition 2.2}
The {\sl degree} of $A\in\Cal B(H)$
is defined by
$$
\deg(A) = \sup_n\ {\text {rank}}(P_nA - AP_n).
$$
\endproclaim

The degree of an operator can
be any nonnegative integer
or $+\infty$, and it is clear
that $\deg(A) = \deg(A^*)$.
Operators of finite degree are
an abstraction
of band-limited matrices.
For example, suppose we have a proper
filtration arising from an
orthonormal basis as in $2.1$.
If the matrix $(a_{ij})$ of an
operator $A\in\Cal B(H)$
relative to this basis satisfies
$$
a_{ij} = 0\quad \text{when}\ |i-j| > k,
$$
then $\deg(A)\leq k$.  The degree function has
a number of natural properties,
the most notable being
$$
\deg(AB) \leq \deg(A) + \deg(B).
$$
The set of all operators of
finite degree is a unital
$*$-subalgebra of $\Cal B(H)$
which is dense in the strong
operator topology \cite{3}.

There is a somewhat larger
$*$-algebra of operators that
one can associate to a filtration
$\Cal F$.  This is a
Banach $*$-algebra relative to
a new norm, and it plays a key
role in the results to follow.
This Banach algebra
is defined as follows.
Suppose that an operator
$A\in\Cal B(H)$ can be decomposed
into an infinite sum of
finite degree operators
$$
A = A_1 + A_2 +\dots\tag{2.3}
$$
where the sequence $A_n$
satisfies the condition
$$
s = \sum_{n=1}^\infty (1+\deg(A_n)^{1/2})\|A_n\|
< \infty.  \tag{2.4}
$$
Notice that the sum indicated
in $2.3$ is absolutely convergent
in the operator norm since $2.4$ implies
$$
\sum_{n=1}^\infty \|A_n\| \leq s <\infty.
$$
We define $|A|_\Cal F$ to be the
infimum of all numbers
$s$ associated with decompositions
 of $A$ of this kind.
If $A$ cannot be so decomposed
then $|A|_\Cal F$ is defined
as $+\infty$.  We define
$$
\Cal D_\Cal F = \{A\in\Cal B(H):
|A|_\Cal F <\infty\}.
$$
$(\Cal D_\Cal F,|\cdot|_\Cal F)$
is a unital Banach algebra
which contains all finite degree
operators, and the operator
adjoint defines an isometric
involution in $\Cal D_\Cal F$.
The two norms currently in view
are related by
$$
\|A\| \leq |A|_\Cal F.
$$

While it is not so easy to calculate
the norm in $\Cal D_\Cal F$,
it is usually very easy to obtain
effective estimates.
For example, let $\Cal F$
be a proper filtration associated
with a unilateral orthonormal
basis as in $2.1$, let $A$ be an
operator on $H$ and let
$(a_{ij})$ be the matrix of $A$
relative to this basis.
Letting
$$
d_k = \sup_i\ |a_{i,i+k}|
$$
denote the supnorm of the $k$th
diagonal, $k=0,\pm 1,\pm 2,\dots$,
then we have
$$
|A|_\Cal F \leq \sum_{k=-\infty}^{+\infty}
(1 + (2|k|)^{1/2})d_k.
$$
See \cite{3}.  Thus, $A$ will belong
to $D_\Cal F$
whenever the diagonals of $(a_{ij})$
die out fast enough so that
$$
\sum_{-\infty}^{+\infty}|k|^{1/2}d_k < \infty.
$$

We now indicate the role of the Banach
algebra $D_\Cal F$
in computing spectra.  Let
$A\in\Cal B(H)$ be a self adjoint
operator and consider the sequence
of matrices $A_n$ obtained
by compressing $A$ along the filtration
$\Cal F = \{H_1\subseteq H_2\subseteq \dots\}$:
$$
A_n = P_n A\restriction_{H_n}.
$$
We are interested in
certain asymptotic quantities
that can be computed (at
least in principle) from the
sequence of finite dimensional
spectra $\sigma(A_1), \sigma(A_2),\dots$.
The simplest one
is the set $\Lambda$, which consists of
all real numbers with the property that
there is a sequence of eigenvalues
$\lambda_n\in\sigma(A_n)$, $n=1,2,\dots$
such that
$$
\lim_{n\to\infty}\lambda_n = \lambda.
$$
$\Lambda$ is in some sense the set of
{\it limit points}
of the sequence of sets $\sigma(A_1),
\sigma(A_2),\dots$.
Notice however that $\Lambda$ is
smaller than the
topological limit
superior of sets because limits
along subsequences of
$\sigma(A_1), \sigma(A_2),\dots$ do
not qualify for
membership in $\Lambda$.
It is easy to see that $\Lambda$ is
closed, and in general
$\Lambda$ contains the spectrum of
 $A$ \cite{3}.  In particular,
$\Lambda$ is never empty.  The points of
$\Lambda$ are classified as follows.

For every subset
$S\subseteq\Bbb R$ and every $n=1,2,\dots$
we write
$N_n(S)$ for the number of
eigenvalues of $A_n$ which belong
to $S$...and of course one counts
multiple eigenvalues according
to their multiplicity.

\proclaim{Definition 2.5} A real
number $\lambda$ is called
a transient point if there is an
open set $U$ containing
$\lambda$ such that
$$
\sup_{n\geq 1}N_n(U) \leq M < \infty.
$$
$\lambda$ is called an essential
point if for every
open set U containing $\lambda$ we have
$$
\lim_{n\to\infty}N_n(U) = \infty.
$$
\endproclaim

\remark{Remarks}
With every transient point
$\lambda$ of $\Lambda$
there is an associated
pair of positive integers $p\leq q$
with the property
that for every sufficiently
small neighborhood $U$ of $\lambda$
one has the behavior
$$
p\leq N_k(U)\leq q
$$
for large enough $k\geq n=n_U$, and moreover
both extreme values $p$ and $q$ are taken on infinitely
many times by the sequence $N_1(U),N_2(U),\dots$.

The set of all
essential points is a subset of
$\Lambda$ which we denote
by $\Lambda_e$.  Again, it can
be seen that $\Lambda_e$
{\it is a closed set which contains
the essential spectrum of}
$A$ \cite{3}.
At this level of generality,
there does not
appear to be much more that
one can say.  For instance, there
are examples which show that
the both inclusions
$\sigma(A)\subseteq \Lambda$
and $\sigma_e(A)\subseteq \Lambda_e$
can be proper (see the appendix of \cite3).
Other examples show that $\Lambda$ can
contain points which are neither
essential nor transient;
for every small neighborhood $U$
of such a point one can find
subsequences $n_1<n_2<\dots$
for which the sequence
of positive integers
$N_{n_1}(U), N_{n_2}(U),\dots$
stays bounded, and other
subsequences $m_1<m_2\dots$ such
that $N_{m_k}(U)\to\infty$
 as $k\to\infty$.  Fortunately,
the following result implies
that in the reasonable cases
we will not find this
kind of instability.
\endremark

We can now state our
 main general result, which shows
how one must choose a
filtration in order to compute
the essential spectrum
of a self adjoint operator.

\proclaim{Theorem 2.6}  Let
$\Cal F$ be a proper filtration and
let $A$ be a self adjoint
operator which belongs
to the Banach $*$-algebra
$\Cal D_\Cal F$.
Then $\Lambda_e$ coincides
with the essential spectrum
of $A$.  Moreover, every
point of $\Lambda$ is either
transient or essential.
\endproclaim

\remark{Remarks}
In the following section we
will encounter
tridiagonal operators $A$
which are defined in terms of
a bilateral orthonormal basis
$\{e_n: n= 0, \pm1, \pm2,\dots\}$
by
$$
Ae_n = e_{n-1} + d_ne_n + e_{n+1},\tag{2.7}
$$
where $d_n$ represents a bounded
sequence of reals.  Theorem
2.6 shows how the essential
spectrum of $A$ can be computed
in terms of the eigenvalue
distributions of the sequence of
$(2n+1)\times(2n+1)$ matrices
obtained by compressing $A$ along the
sequence of subspaces $H_n=[e_{-n},e_{-n+1},\dots,e_{n-1},e_n]$,
$n=1,2,\dots$.  While Theorem
2.6 says nothing about {\it rates}
of convergence, our experience
with the operators of
section 4 has shown that
convergence is rapid; in fact, it
is fast enough to allow the
construction of excellent
pictures of $\sigma_e(A)$
on desktop
computers in a minute or
two \cite{4}.  Nevertheless, there
remains an important problem
of obtaining an appropriate
definition of the ``rate" of
convergence of such approximations,
and the estimation of this rate.
The paper \cite{3}
does not address these issues
of error definition and
estimation.
\endremark

In carrying out computations,
it is usually more convenient
to deal not with the above
``bilateral" sections $A_n$
but with smaller ``unilateral"
sections.  More precisely,
these are defined in terms of an
{\it improper} filtration
$\{H_1^+\subseteq H_2^+\subseteq\dots\}$,
$H_n^+=[e_1,e_2,\dots,e_n]$,
and the corresponding
compressions of the operator $A$
$$
A_n^+=P_n^+A\restriction_{H_n^+},
$$
$P_n^+$ denoting the projection
onto $H_n^+$.  Since the
filtration $\{H_n^+\}$ is
improper, Theorem 2.6 does not
allow one to draw conclusions
about the essential
spectrum of $A$, but rather the
essential spectrum of
the operator $A^+$ obtained by compressing $A$ to the
subspace $H^+$ generated by
$\cup_nH_n^+$.
Thus there remains a significant
problem of relating the essential
spectrum of $A$ to
that of $A^+$.  These issues will
be taken up in
section 4 below.

\subheading{3. Discretized Schr\"odinger operators}
In this section we will indicate
how a one-dimensional
quantum system should be discretized
in order to carry
out numerical computations.
We will find that the
``discretized" canonical commutation
relations generate
a \cstar\ which is isomorphic to one of the
noncommutative spheres of Bratteli,
Elliott, Evans and
Kishimoto (with parameter
related to the numerical step
size), and we will find that
the resulting discretized
Hamiltonian is a bounded
self-adjoint operator which
is amenable to the methods
of the preceding
section.

Most one-dimensional quantum
systems are modelled on
the Hilbert space $L^2(\Bbb R)$.
The canonical operators
are the unbounded self-adjoint
operators $P, Q$ defined
on appropriate domains in $L^2(\Bbb R)$ by
$$
\align
P &= \frac{1}{i}\frac{d}{dx},  \tag{3.1}\\
Q &= \text{Multiplication by }x.
\endalign
$$
They obey the canonical commutation relations on an
appropriate common domain
$$
PQ - QP = \frac{1}{i}{\bold 1}.
$$
The time development of the
system is described by a
one-parameter unitary group of the form
$$
W_t = e^{itH},\qquad t\in\Bbb R,
$$
where $H$ is the
Hamiltonian of the system
$$
H = \frac{1}{2}P^2 + \phi(Q),\tag{3.2}
$$
$\phi$ being a real-valued
continuous function of a real
variable which represents the
 potential of the classical
system being quantized.

In order to carry out numerical
computations one first
has to replace the differential operator $H$ with an
appropriately discretized version of itself.
Moreover, one has to decide
how this should be done so
as to conform with the basic
principles of numerical
analysis while at the same
time preserving the essential
features of quantum mechanics
(i.e., the uncertainty
principle).  In \cite1
we presented arguments which
we believe justify the
following procedure.

One first settles on a numerical
step size $\sigma$.  This
can be regarded as a small
positive rational
number, whose size
represents the smallest
time increment to be used in
the difference equations
that replace differential
equations.  One then
discretizes {\it both\/} operators
$P$ and $Q$, and finally
uses the formula 3.2 to define
the corresponding discretized
version of the Hamiltonian.

In more detail, we replace the
differential
operator $P$ with the bounded
self-adjoint difference operator
$P_\sigma$ defined by
$$
P_\sigma f(x) = \frac{f(x+\sigma)-f(x-\sigma)}{2i\sigma},
\qquad x\in\Bbb R.
$$
Noting that the one-parameter group of translations
$V_tf(x)=f(x+t)$ is generated
by $P$ in the sense that
$V_t = e^{itP}, t\in\Bbb R$, we have
$$
P_\sigma = \frac{1}{2i\sigma}(e^{i\sigma P}-e^{-i\sigma P}) =
\frac{1}{\sigma}\sin(\sigma P).
$$
Now we must discretize $Q$ but
we must be careful to
do it in a way that preserves
the uncertainty principle
insofar as that is possible.
In section 3 of \cite1, we argued
that this requirement
imposes a very strong restriction
on the possible choices
for the discretized $Q$, and that
in fact the only
``correct" choice is given by
$$
Q_\sigma = \frac{1}{2i\sigma}(e^{i\sigma Q}-e^{-i\sigma Q})
=\frac{1}{\sigma}\sin(\sigma Q).
$$
More explicitly, we have
$$
Q_\sigma f(x) = \frac{1}{\sigma}\sin(\sigma x)f(x),
\qquad x\in\Bbb R, f\in L^2(\Bbb R).
$$
The resulting discretized Hamiltonian
is then defined as follows:
$$
H_\sigma = \frac{1}{2}P_\sigma^2 + V(Q_\sigma).  \tag{3.3}
$$

Once we have the operator $H_\sigma$ we are in position
to carry out numerical computations with
the quantum system.
For example, if the system is described
at time $t$ by a
wave function $f\in L^2(\Bbb R)$
then the state of the
system at time $t+\Delta t$ is
approximated by the wave
function
$$
g = ({\bold 1}+i\Delta t H_\sigma)f =
f + i\Delta t H_\sigma f.
$$
We caution the reader that, while the
preceding formula is convenient
for illustrating one way to make use
of the discretized Hamiltonian, it
must be modified appropriately in
order to correctly model
the dymanical group in practice
because the operator $f\mapsto g$ is not
unitary.  Readers interested in
carrying out numerical computations
can find a discussion of closely
related issues in
\cite{14, pp 662--663}.

Let $\Cal D_\sigma$ be the \cstar\
generated by the set
of operators $\{P_\sigma, Q_\sigma\}$.
$\Cal D_\sigma$
is the norm-closed linear span of all
finite products
of terms involving either $P_\sigma$
or $Q_\sigma$.  It
is not obvious that $\Cal D_\sigma$
contains the identity
operator but that is the case.
$\Cal D_\sigma$ is
the discretized counterpart
of the algebra of observables,
and notice that it contains
the operator $H_\sigma$.
Thus it is important to understand
the structure of $\Cal D_\sigma$.

In fact, $\Cal D_\sigma$ is a
simple unital \cstar\ which is
isomorphic to one of the

non-commutative spheres of Bratteli
et al \cite{9, 10, 11}.  Moreover,
while the operators
$P_\sigma, Q_\sigma$ no longer
satisfy the canonical
commutation relations, they do
obey a more subtle discretized
form of the CCRs, and more generally
it is the universal
\cstar\ associated with these
``discretized" CCRs that is
naturally associated with the
non-commutative spheres.
The reader is referred to \cite{2}
for a detailed discussion
of these and related issues.

In particular, if one is interested
in computing the spectrum
of $H_\sigma$ then one is free
to choose any convenient representation
of $\Cal D_\sigma$ and compute
the spectrum of $H_\sigma$ in
that representation.  Since
$\Cal D_\sigma$ is simple,
the spectrum does not depend
on the representation chosen.
Actually, the most convenient
realization of $H_\sigma$ is
one in which it is a tridiagonal
operator.  In this case,
it is more appropriate to work
with a subalgebra of
$\Cal D_\sigma$ which contains
$H_\sigma$.  The precise
statement follows.

\proclaim{Proposition 3.4}
Let $\Cal A$ be the
\cstar\ generated by
$P_\sigma^2$ and $Q_\sigma$, and
let $K$ be a Hilbert space
spanned by a bilateral
orthonormal set $\{e_n: n\in\Bbb Z\}$.
Then there is
a faithful representation
$\pi:\Cal A\to\Cal B(K)$
such that $\pi(H_\sigma)$
 has the form
$$
\pi(H_\sigma) = aT + b{\bold 1}\tag{3.5}
$$
where $a = 1/{8\sigma^2}$,
$b=-1/{4\sigma^2}$, and T
is the tridiagonal operator
$$
Te_n =
e_{n-1} + 8\sigma^2\phi(\frac{1}{\sigma}\sin(2n\sigma))e_n +
e_{n+1},\tag{3.6}
$$
$n=0,\pm1, \pm2,\dots$.
\endproclaim
\demo{proof}
Consider the unitary operators $U$,
$V$ defined by
$$
\align
Uf(x) &= e^{i\sigma x}f(x)\\
Vf(x) &= f(x+\sigma).
\endalign
$$
Then
$$
\align
P_\sigma &= \frac{1}{2i\sigma}(V - V^{-1}),\\
Q_\sigma &= \frac{1}{2i\sigma}(U - U^{-1}),
\endalign
$$
and notice that
$$
P_\sigma^2 = -\frac{1}{4\sigma^2}
(V^2 + V^{-2}) +\frac{1}{2\sigma^2}.
$$
Let $\Cal B$ be the \cstar\ generated
by $U$ and $V^2$.  Clearly
$\Cal A$ is contained in $\Cal B$, and
because of the commutation
relation
$$
V^2 U = e^{2i\sigma^2}UV^2
$$
and the fact that $\sigma$ is a
positive rational number,
it follows that $\Cal B$ is an
irrational rotation \cstar.
We will define a representation
$\pi_1$ of $\Cal B$ on
$K$; the required representation
of $\Cal A$ is obtained by
restriction.

In order to specify a representation
of $\Cal B$ on $K$, it
is sufficient to specify a pair of
unitary operators $S$, $D$
on $K$ satisfying
$$
SD = e^{2i\sigma^2}D S.  \tag{3.7}
$$
$\pi_1$ is then uniquely
defined by specifying that
$\pi_1(V^2) = S$, $\pi_1(U)=D$.  Let
$$
\align
Se_n &= -e_{n-1},\\
De_n &= e^{2in\sigma^2}e_n
\endalign
$$
$n=0,\pm1, \pm2,\dots$.
Clearly $S$ and $D$ are
unitary operators and the
reader can verify 3.7 directly.  Hence
there is a representation $\pi_1$
of $\Cal B$ with
the stated properties.  Noting that
$$
\align
\pi_1(P_\sigma^2) &=
\pi_1(-\frac{1}{4\sigma^2}(V^2+V^{-2}) +\frac{1}{2\sigma^2}{\bold 1})\\
&= -\frac{1}{4\sigma^2}S -
\frac{1}{4\sigma^2}S^{-1} + \frac{1}{2\sigma^2}{\bold 1}
\endalign
$$
and that $\pi_1(Q_\sigma)$ is the diagonal operator
$$
\pi_1(Q_\sigma)e_n =
\frac{1}{2i\sigma}(e^{2in\sigma^2} - e^{-2in\sigma^2})e_n
= \frac{1}{\sigma}\sin(2n\sigma^2)e_n,
$$
we may conclude that
$$
\pi_1(\phi(Q_\sigma))e_n = \phi(\frac{1}{\sigma}\sin(2n\sigma^2))e_n
$$
for every $n=0,\pm1,\pm2,\dots$.
Combining these two formulas, we
find that the image of $H_\sigma$ is given by
$$
\pi_1(H_\sigma) = \frac{1}{2}\pi_1(P_\sigma^2) +
\pi_1(\phi(Q_\sigma)),
$$
which has the form spelled out in 3.5 and 3.6\qed
\enddemo

\remark{Remark 3.8}  Consider
the diagonal sequence
$$
d_n = 8\sigma^2\phi(\frac{1}{\sigma}\sin(2n\sigma^2)),
\qquad n\in\Bbb Z
$$
appearing in the formula 3.6.
We want to point out that if
the function $\phi$ is
continuous and not constant on the
interval $[-\frac{1}{\sigma},+\frac{1}{\sigma}]$ and if
$\sigma^2$ is not a rational
multiple of $\pi$, then the
sequence $(d_n)$ is
{\it almost periodic} but {\it not
periodic}.  Indeed, the sequence
$n\mapsto \sin(2n\sigma^2)$ is
almost periodic because it is a
linear combination of
complex exponentials of the form
$e^{i\alpha n}$ where $\alpha$
is a real number.  Since the set
of all almost periodic sequences
form a commutative \cstar\ it is closed
under the continuous functional
calculus, and therefore the sequence
$(d_n)$ must be almost periodic.
Note too that $(d_n)$ cannot be periodic.
For if there did exist integers
$p\geq 1$ and $n$ such that
$$
d_{n+kp} = d_n, \qquad k=0,\pm1, \pm2,\dots
$$
then since $\sigma^2$ is not a
rational multiple of $\pi$ the numbers
$$
\sin(2n\sigma^2 + 2kp\sigma^2),\qquad k\in\Bbb Z
$$
would fill out a dense set in the
interval $[-1,+1]$, and hence
$\phi$ would have to be
constant on the interval
$[-\frac{1}{\sigma},+\frac{1}{\sigma}]$.

The material presented in this section
leads toward a significant
conclusion.  {\it The problem of
computing the spectrum of the discretized
Hamiltonian of a one dimensional
quantum system can be reduced
to the problem of computing the
spectrum of a self-adjoint tridiagonal
operator of the form}
$$
Te_n = e_{n-1} + d_ne_n + e_{n+1},
\qquad n\in\Bbb Z\tag{3.9}
$$
{\it where} $\{d_n: n\in\Bbb Z\}$
{\it is a bounded almost periodic sequence
of reals which is not periodic}.
For such operators one can
work with either unilateral sections
or bilateral sections, as we
will see in the following section.
For example, if for $n\geq 1$ we let
$T_n$ be the compression of $T$ onto
the linear span of $\{e_1,e_2,\dots,e_n\}$
then even though we are in effect working
with an improper filtration, one may
apply a suitable variation of
Theorem 2.6 to conclude that the spectrum
of $T$ is the set of all essential
points associated with the
sequence of self-adjoint matrices
$T_1, T_2, \dots$.
\endremark

In the following section, we will
obtain more precise information
about the distribution of the
eigenvalues of the sequence $T_n$,
$n=1,2,\dots$.

\subheading{4. Limits, simple
\cstar s, and traces}
Let $(d_n)_{n\in\Bbb Z}$ be a
bounded almost periodic
sequence of reals which is
{\it not} periodic and let
$T$ be the tridiagonal
operator of 3.9
$$
Te_n = e_{n-1} + d_ne_n + e_{n+1},
\qquad n\in\Bbb Z,\tag{4.1}
$$
$\{e_n: n\in \Bbb Z\}$ being an
orthonormal basis for
a Hilbert space $H$.  We will
show that the eigenvalue
distributions of the $n\times n$
sections of $T$
actually converge (in the weak$^*$
topology of
measures on the real line) to a
probability measure
$\mu_T$.  Thus it becomes important
to understand
the nature of this limiting measure
$\mu_T$.  We show
that $\mu_T$ is associated
with a tracial state on
a certain \cstar\ $\Cal A_T$
associated with $T$.
$\Cal A_T$ certainly contains
$T$ but it is much
larger than the \cstar\ generated
by $T$; indeed,
$\Cal A_T$ is a simple \cstar\
having a {\it unique}
tracial state.

Let $S$ be the bilateral shift
defined on $H$ by
$$
Se_n=e_{n+1} \qquad n\in\Bbb Z
$$
and let $D$ be the diagonal operator associated
with the sequence $(d_n)$,
$$
De_n = d_ne_n,\qquad n\in \Bbb Z.
$$
$\Cal A_T$ is defined as the \cstar\
generated by
both operators $D$ and $S$.
$\Cal A_T$ is a separable
unital \cstar\ which contains $T$,
and of course
$\Cal A_T$ depends on the particular
choice of
diagonal sequence $(d_n)$.
When $T$ is a
discretized Hamiltonian as in
the previous section,
$\Cal A_T$ will depend on both
the numerical step
size $\sigma$ and the potential
$\phi$.  Nevertheless,
in all cases we have

\proclaim{Theorem 4.2}$\Cal A_T$
is a simple unital
\cstar -subalgebra of
$\Cal B(H)$ which has a
unique tracial state.
\endproclaim

\remark{Remark}
A {\it tracial state} of
$\Cal A_T$ is a linear functional
$\tau:\Cal A_T \to \Bbb C$
satisfying
$$
\align
&\tau(X^*X)\geq 0, \quad{\text and}\\
&\tau(XY) = \tau(YX)
\endalign
$$ for every $X,Y\in\Cal A_T$, and which
is normalized so that $\tau({\bold 1})=1$.
Theorem 4.2
is proved in proposition 3.2 of \cite{5}.
\endremark

More generally, suppose we are
given an arbitrary concrete
\cstar\ $\Cal A\subseteq\Cal B(H)$.
We need
to single out the filtrations
that are ``compatible"
with $\Cal A$.  Let
$\Cal F=\{H_1\subseteq H_2\subseteq\dots\}$
be a filtration of $H$ which may be
improper, and
put
$$
H_+ = \overline{\cup_nH_n}.
$$
As in section 2 we may speak of the degree of an
operator $X\in \Cal B(H)$ relative to
the filtration $\Cal F$.
Following \cite{3}, we
say that $\Cal F$ is an $\Cal A$-filtration if
the set of finite degree operators
which belong
to $\Cal A$ is norm-dense in $\Cal A$.

In connection with the operators
$T$ of 4.1, we
will consider the filtration
$\Cal F_T = \{H_n\}$
where $H_n=[e_1,e_2,\dots,e_n]$.
This is an
improper filtration for which
$$
H_+ = [e_1,e_2,\dots].
$$
More significantly, we have

\proclaim{Proposition 4.3}
$\Cal F_T$ is an $\Cal A_T$-filtration.
\endproclaim

4.3 is a simple consequence
of the fact that
the finite degree operators in
$\Cal A_T$ form a 
of $\Cal A_T$, and that the operators
$D$ and $S$ have respective degrees
$0$ and $1$ (see the proof of
Theorem 3.4 of \cite{5}).
For every $n\geq 1$ let $T_n$
be the compression of
$T$ to $H_n$.  Relative to the
obvious basis for $H_n$,
the matrix of $T_n$ is
$$
\pmatrix
d_{1}&1&0&\hdots&0&0\\
1&d_{2}&1&\hdots&0&0\\
0&1&d_{3}&\hdots&0&0\\
\vdots&\vdots&\vdots&\ddots&\vdots&\vdots\\
0&0&0&\hdots&d_{n-1}&1\\
0&0&0&\hdots&1&d_n
\endpmatrix.
$$

\noindent There are two issues
that need to be understood.
The first has to do with the
relation between operators
in $\Cal A_T$ and their compressions
to $H_+$; the
second requires relating the trace
on $\Cal A_T$ to
the limit of the eigenvalue
distributions of the
sequence of matrices $T_1, T_2,\dots$.

In order to discuss the first
of these, we consider
a more general setting in which
we are given a unital
\cstar\ $\Cal A\subseteq\Cal B(H)$ and an improper
filtration
$\{H_1\subseteq H_2\subseteq\dots\}$.
We will consider the space
$\Cal A_+$ of all
operators on $H_+$ having the form
$$
P_+A\restriction_{H_+} + K,
$$
where $A\in\Cal A$ and $K$
is a compact operator
on $H_+$, $P_+$ denoting
the projection on $H_+$.
We will write $\Cal K_+$
for the algebra of all
compact operators on $H_+$.

\proclaim{Theorem 4.4}
Assume that $\Cal A$ has
a unique tracial state
$\tau$, let $\{H_1\subseteq H_2\subseteq\dots\}$
be an $\Cal A$-filtration
and assume that $H_+$
has the following property
$$
A\restriction_{H_+} = {\text compact}
\implies A = {\text compact}
\tag{4.5}
$$
for every operator $A\in\Cal A$.
Then $\Cal A_+$
is a \cstar\ containing
$\Cal K_+$ which has a
unique tracial state $\tau_+$.
$\tau_+$ is related
to $\tau$ by way of
$$
\tau_+(P_+A\restriction_{H_+} + K)
= \tau(A),
\qquad A\in\Cal A, K\in \Cal K_+.
$$
Moreover, the natural map of
$\Cal A$ to the quotient
$\Cal A_+/\Cal K_+$ given by
$$
A\mapsto P_+A\restriction_{H_+} + \Cal K_+
$$
is an isomorphism of \cstar s:
$$
\Cal A\cong \Cal A_+/\Cal K_+.
$$
\endproclaim

\remark{Remarks}
The argument required here can be found
in the proof of Theorem 2.3
of \cite{5}.  These results depend on the following
relationship that exists between the operators in
$\Cal A$ and the projection $P_+$ associated with
an $\Cal A$-filtration:
$$
A\in\Cal A\implies P_+A - AP_+ \in \Cal K
$$
see \cite {5, proposition 2.1}.

Theorem 4.4 implies that we have
a short exact sequence
of \cstar s
$$
0\to\Cal K_+\to \Cal A_+\to \Cal A\to 0.
$$
This shows that the structure of
$\Cal A_+$ is somewhat
analogous to the structure of
the Toeplitz \cstar\ $\Cal T$,
$$
0\to\Cal K\to \Cal T\to C(\Bbb T)\to 0,
$$
except that in our applications
the quotient \cstar s
$\Cal A_T$ are highly noncommutative.

Finally, we remark that the
hypothesis 4.5 is satisfied
for our examples $\Cal A_T$ because
$\Cal A_T$ is a simple
\cstar\ with unit and
$H_+$ is infinite dimensional
(see the proof of Theorem
3.4 of \cite{5}).
\endremark

Now we can define the measure
$\mu_T$ alluded to at
the beginning of this section.
Let $\tau$ be the
unique tracial state of $\Cal A_T$.
Using the
functional calculus for bounded
self-adjoint operators,
we can define a positive
linear functional on
$C_0(\Bbb R)$ by
$$
f\in C_0(\Bbb R)\mapsto \tau(f(T)).
$$
By the Riesz-Markov theorem,
there is a unique positive
measure $\mu_T$ on
$\Bbb R$ such that
$$
\int_{-\infty}^{+\infty}f(x)\,d\mu_T(x) =
\tau(f(T)),\qquad f\in C_0(\Bbb R).
$$
$\mu_T$ is obviously a
probability measure
whose support is contained
in the spectrum of $T$.
$\mu_T$ is called the {\it
spectral distribution} of
$T$.  Because $\Cal A_T$ is
simple $\tau$ must be
a faithful trace, and hence
the closed support of
$\mu_T$ is {\it exactly\/}
the spectrum of $T$.

The preceding discussion,
together with the general results
of \cite{3, section 4} can
be applied to obtain the following
result, which is
Theorem 3.4 of \cite{5}.

\proclaim{Theorem 4.6}
For every positive integer $n$
let $\lambda_1^n < \lambda_2^n <\dots <\lambda_n^n$
be the eigenvalue list of the
symmetric $n\times n$ matrix
$$
\pmatrix
d_{1}&1&0&\hdots&0&0\\
1&d_{2}&1&\hdots&0&0\\
0&1&d_{3}&\hdots&0&0\\
\vdots&\vdots&\vdots&\ddots&\vdots&\vdots\\
0&0&0&\hdots&d_{n-1}&1\\
0&0&0&\hdots&1&d_n
\endpmatrix.
\tag{4.7}
$$
Then for every $f\in C_0(\Bbb R)$ we have
$$
\lim_{n\to\infty}\frac{1}{n}(f(\lambda_1^n)+
f(\lambda_2^n)+\dots,
+f(\lambda_n^n)) =
\int_{-\infty}^{+\infty} f(x)\,d\mu_T(x).
$$
\endproclaim

\remark{Remarks}
Theorems 4.6 and 2.6 together
provide rather precise
information about the
rate at which the eigenvalues
of the matrices 4.7
accumulate at points in and
out of the spectrum of $T$.  For
example, let
$\lambda\in\sigma(T)$ and
let $I$ be
an open interval
containing $\lambda$.  Then
$\mu_T(I)>0$, and if
$\alpha$ and $\beta$ are
chosen close to $\mu_T(I)$ in such a way that
$\alpha < \mu_T(I) < \beta$, then
the number $N_n(I)$
of eigenvalues of $T_n$ which
belong to $I$ will
satisfy the inequalities
$$
\alpha n\leq N_n(I) \leq \beta n
$$
for all sufficiently large $n$.

For the operators $T$ of
4.1, it is not hard
to show that the essential
spectrum of $T$ is
identical with $\sigma(T)$.
Thus we can apply
Theorem 2.6 above to conclude that
if $\lambda$ does {\it not}
belong to
$\sigma(T)$ then for every
sufficiently small
open interval containing
$\lambda$, the sequence of
numbers $N_1(I), N_2(I),\dots$
actually
remains bounded.  This
behavior is quite visible
in the pictures generated by
the program \cite{4}.

Some additional
remarks relating to computational
issues can be found
in the concluding paragraphs of \cite{5}.
\endremark

%
%

\enddocument